\documentclass{elsart}
\journal{Chaos, Solitons \& Fractals}
\volume{25/1, pp. 29-32}
\usepackage{graphicx,amssymb}

\begin{document}

\begin{frontmatter}

\title{Some non-conventional ideas about algorithmic complexity}
\author{Germano D'Abramo}
\address{Istituto di Astrofisica Spaziale e Fisica Cosmica,\\ 
Area di Ricerca CNR Tor Vergata, Roma, Italy.\\
E-mail: Germano.Dabramo@rm.iasf.cnr.it}

\date{(Accepted)}

\begin{abstract} 

In this paper the author presents some non-conventional thoughts on the 
complexity of the Universe and the algorithmic reproducibility of the 
human brain, essentially sparked off by the notion of algorithmic 
complexity. We must warn that though they evoke suggestive scenarios, they 
are still quite speculative.

\end{abstract}

\end{frontmatter}

\section{Introductory remarks}

The algorithmic (program-size, or Kolmogorov) complexity of a binary
string $s$ is defined as the size in bits of the smallest computer program
able to generate it:
$$H(s)\equiv \min_{U(p)=s} |p|$$
where $p$ is a program string used by a universal computer $U$ to produce 
the sequence $s$ (Chaitin~\cite{sa,ait}).

Let us now consider the following algorithm of size $\lfloor\log_2 
N\rfloor +k+1$ (see also the Appendix). It lists in order of their size 
all strings of length less than or equal to $N$ bits (there are 
$2^{N+1}-2$ of them\footnote{ That is, we sum up the number of all 
possible strings of one bit ($2$), that of all possible strings of two 
bits ($2^2$), and so on up to $N$ bits: $$\sum_{i=1}^{N} 2^i = 
2+2^2+2^3+...+2^{N}=2^{N+1}-2.$$}), writing them one after the other on a 
computer file (or on a sheet of paper). The output should look like

\begin{displaymath}
0\,\,1\,\,00\,\,01\,\,10\,\,11\,\,000\,\,001\,\,010\,\, 
011\,\,100\,\,101\,\,110\,\,111\,\,...
\end{displaymath}

Let us assume also that $N\gg\lfloor\log_2 N\rfloor +k+1$. Thus, since
$$\frac{2^{\lfloor\log_2 N\rfloor +k+2}-2}{2^{N+1}-2}\ll 1,$$ where
$2^{\lfloor\log_2 N\rfloor +k+2}-2$ is the total number of different
program strings which could be obtained with a number of bits less than or
equal to $\lfloor\log_2 N\rfloor +k+1$, then among all the produced
strings there is surely at least one that is more complex than our
algorithm, i.e.~it can't be produced by any program of size less than or
equal to $\lfloor\log_2 N\rfloor +k+1$ bits, by definition of algorithmic
complexity. Suddenly a paradox appears: an algorithm of size equal to
$\lfloor\log_2 N\rfloor +k+1$ bits is able to write a list of strings
which contains at least one that is more complex than $\lfloor\log_2
N\rfloor +k+1$ bits, namely than the algorithm itself.

It is even more striking to consider a similar but simpler algorithm (the
same as before but without the {\tt if} condition, and thus having
constant size, nearly equal to $k$ bits; more in the Appendix) which lists
every natural number in binary notation, endlessly. In this context, we
have that every binary string, of any length, is generated by this almost
trivial algorithm. And, although every single string could be extremely
complex (for instance, like that coding the collected works of Giacomo
Leopardi) the whole, infinite set has a ridiculous algorithmic complexity.

Think for a moment to a god less `complex' than a small portion of what it
created! Obviously, we are not claiming that our Universe is {\it tout
court} identical to a computer algorithm; what we want to suggest is that
the concept of something simpler generating something incredibly more
complex is not completely weird and unfounded, where simplicity and
complexity must be intended in terms of algorithmic complexity.

Think for a while to the mocking and frustrating possibility that this
applies to our real world too. Mankind is searching for the ultimate
meaning of things, expecting to find it through a wider and more complex
description than we currently have. But, maybe, that assumption could be
simply wrong, and the ultimate reason, the origin of complex things might
be in a simpler thing, as it happens with our listing algorithm. Saying it
in other words, complexity might be simply a by-product and any wider and
more complex description, far from being a step toward the final
explanation of things, might exist with no finality, might be simply
self-aimed. Maybe, Who/What created the Universe is not Omniscient and
Almighty like we suppose It should be, but, let us say, It might be
simpler than bacteria. Tegmark~\cite{te} also suggested the possibility of
a Universe with almost no information content if taken as a whole, and he
did it from the point of view of quantum mechanics, while
El~Naschie~\cite{el}, using the Newton non-dimensional gravity constant
$\overline{\alpha}_G$ as a measure of complexity for the Universe, found
that the information dimension of the Universe is 128, nearly equal to the
inverse of the Sommerfeld electromagnetic fine structure constant measure
at the electro-weak scale $\overline{\alpha}_{ew}$.

Besides, the property of being explainable to mankind might be distributed
by chance over the things of our world. Some strings/patterns are
reproducible by shorter (therefore simpler and more familiar) programs,
and thus they are `explainable', while others are not and they are
perceived as random by us (see, for example, the model of inductive
inference by Solomonoff~\cite{so}; and the notion of random string,
Chaitin~\cite{sa}). But there might be no design under the distribution of
what is explainable and what is not. Likewise, some strings of our listing
program are reproducible by shorter program, some other not, that's all!

A similar argument on the (un)explainability of the reality was proposed
by C.~S.~Calude in terms of {\em lexicons} (see Calude and
Meyerstein~\cite{cm} and references therein). A {\em lexicon} is the
infinite expansion of some real number (e.g. the infinite binary expansion
of 0s and 1s obtained through the tossing of a fair coin) with the
base-independent property of containing every finite string as a
sub-string, infinitely many times. Calude and Meyerstein~\cite{cm} suggest
that the Universe might behave like a lexicon and that we maybe `live' on
a very long finite sequence that is ordered, i.e.~it may be explained
through science, at least partially. But, there is no guarantee that it is
anywhere, anytime the same; the order may suddenly switch to pure
randomness in other portions of the lexicon/Universe.

\section{Probability and the brain}

\begin{quote}

{\em The probability that your brain, specifically its function of giving
rise to your mind, is reproducible by a finite algorithm (e.g.~$N$-bit
long) is arbitrarily close to 0.}

\end{quote}

\noindent For the sake of thought experiment, let us try to provide some
arguments in favor of the above statement. Let us suppose that the
processes of your brain which give rise to your mind are exactly
assimilable to and reproducible by an algorithm of $N$ bits.  Thus, we are
able to algorithmically reproduce your brain (and thus its product, your
mind) as a program running on an `hardware' different from your body in
such a way that it does not suffer from the main limitation of the human
body, namely its relatively short life. In this way, such `brain' could in
principle work for an unlimited time span.

Moreover, let us imagine that as a part of the program simulating your
brain there is a very simple (from the algorithmic point of view) and fast
counter able to enumerate (but not to store) all binary numbers of $c$
bits (there are $2^{c}$ of them) in increasing lexicographical order.

Now, the reproduced mind could {\it think of} a specific decimal number
greater than the size (in bits) of its own generating algorithm, say
$N+k$, and make the counter start counting in increasing order all
possible binary numbers/strings of $N+k$ bits (there are $2^{N+k}$ of
them).

At its own will, the simulated mind can then stop the counter whenever it
wants and print the last enumerated number. The counter could be provided
with a sort of counting completeness indicator, which would give the
percentage of the whole count reached till that moment. This can be of
some help to the simulated mind in choosing when to stop the counter (not
too early for instance, since the first binary strings are surely not very
complex algorithmically).

Remember that the printed number is a binary string of $N+k$ bits, less
than $2^{N+k}$ in size, and it results to be somewhat blindly chosen by
the simulated mind. Therefore, with probability nearly equal to

$$1-\frac{2^{N+1}-2}{2^{N+k}}\sim 1-2^{1-k},$$ 
where $2^{N+1}-2$ is the total number of different strings/programs which
could be obtained with a number of bits less than or equal to $N$, the
simulated $N$-bit long algorithmic brain/mind would be able to generate a
sequence of a complexity greater than $N$ bits.

For $k\gg1$, such probability becomes arbitrarily close to $1$. To recap,
within the hypothesis that your brain (and thus its product, your mind) is
reproducible by an $N$-bit long algorithm, the probability that such
algorithm would be able to generate a string of a complexity greater than
$N$ bits, and thus leading to a logical paradox according to the
definition of algorithmic complexity, is arbitrarily close to $1$.  
Therefore, the probability that the hypothesis will be violated is
arbitrarily close to $1$.

Of course, there is a number of possible critiques to the above argument.  
Someone may argue that for big values of $N$ the described enumeration
(involving $2^{N+k}$ binary numbers of $N+k$ bits) is not physically
feasible (it would require an incredibly huge amount of time), making our
point physically unsound.

Others may argue that the above argument does not eliminate at all the
possibility that our minds are `algorithmic': we might be similar to
machines, behaving predictably like machines (and thus `choosing', through
the procedure described above, a string of consistent algorithmic
complexity), but simply and wrongly believing we are not.

But, even in such cases our argument should be of some interest: though
probably physically unfeasible, our point seems to be in principle
logically and mathematically sound (after all, many trusted mathematical
demonstrations are physically unverifiable, for they involve the concept
of infinity for instance) and maybe it might be an example of a physical
status (i.e.~our brains/minds actually like algorithms) for which we are
able to provide a logically and mathematically sound argument of the
contrary.

\section{Digression}

Let us consider the following device. A mechanical tool reads a decimal
number $N$ in input and tosses an idealized, fair coin $N$ times.
Whenever the result of the toss is a head, such device prints a $1$ on a
long tape, otherwise it prints a $0$. One might think that the {\it
algorithmic size} of this device is proportional to $\log_2 N$ (the size
of the binary expansion of the decimal number $N$), and therefore that it
could be able to generate a sequence of a complexity greater than $\log_2
N$ bits, with probability arbitrarily close to $1$.

In that case, however, the algorithmic size of the device is comparable to
the complexity of the entire physical process, of all the physical laws
(plus all the relevant initial conditions) which make the toss to result
each time in a head rather than a tail, or vice versa. Hence, this device
is not an algorithm {\it embedded} in an electronic computer;  rather, it
operates under the influences of the physical world.  Maybe the same
argument might apply to human mind too: the peculiarity of the brain, in
giving rise to mind and consciousness, might partly originate from its
complex and continuous physical interaction with the surrounding physical
world.

\section*{Acknowledgements}

I am very grateful to Professor~M.~S.~El~Naschie for helpful suggestions 
and comments on the subject of this paper.

\newpage

\section*{Appendix}

A possible (FORTRAN-like) computer code for our listing algorithm: 

\begin{verbatim}

    program list
\end{verbatim}

\noindent \hspace{0.7cm} {\tt N = Number of Bits} (decimal number)

\begin{verbatim}
    s = 0

 1  write(*,*) bin(s)
    if(#Bit(bin(s)).eq.N) goto 2
    s = s+1
    goto 1
 2  stop

    end
  
\end{verbatim}

The size of this algorithm is a constant $k$ (which results from the
coding of its instructions different from the decimal number $N$) plus the
number of bits of the binary coding of $N$, that is $\lfloor\log_2
N\rfloor +1$ (where $\lfloor x\rfloor$ is the floor function of $x$).

A simpler, yet more powerful, listing algorithm would be the following:

\begin{verbatim}

    program list2

    s = 0

 1  write(*,*) bin(s)
    s = s+1
    goto 1

    end
  
\end{verbatim}

Its size is nearly equal to $k$ but it is able to list an infinite number
of binary strings.

\end{document}